\documentclass{pspum-l}

\newtheorem{theorem}{Theorem}[section]
\newtheorem{lemma}[theorem]{Lemma}
\newtheorem{proposition}[theorem]{Proposition}
\newtheorem{corollary}[theorem]{Corollary}

\theoremstyle{definition}

\theoremstyle{remark}
\newtheorem{remark}[theorem]{Remark}

\numberwithin{equation}{section}

\usepackage{hyperref,amsmath,mathabx,stackrel}
\usepackage[all]{xy}
\usepackage{yfonts}
\usepackage{mathrsfs}

\usepackage{tikz-cd}
\usetikzlibrary{arrows.meta}
\usepackage{adjustbox}

\newtheorem*{definition*}{Definition}

\newcommand{\opcit}{{\it op.cit.\/}\ }
\newcommand{\ie}{{\it i.e.\/}\ }

\def\A{{\mathbb A}}
\def\C{{\mathbb C}}
\def\F{{\mathbb F}}
\def\Q{{\mathbb Q}}
\def\R{{\mathbb R}}
\def\Z{{\mathbb Z}}

\def\cA{{\mathcal A}}

\def\cE{{\mathcal E}}
\def\cF{{\mathcal F}}
\def\cG{{\mathcal G}}
\def\cL{{\mathcal L}}
\def\cO{{\mathcal O}}
\def\cH{{\mathcal H}}
\def\cR{{\mathcal R}}
\def\cS{{\mathcal S}}
\def\cT{{\mathcal T}}

\def\tr{{\rm Tr}}
\def\rep{\vartheta}

\def\Mod{{\rm Mod}}

\def\fourier{\F}

\def\sr0{{\cS^{\rm ev}_0}}
\def\sarr{{\cS(\A_\Q)}}
\def\hatz{{\hat\Z^\times}}
\def\rnt{{[0,\infty)\rtimes{\N^{\times}}}}

\def\GL{{\rm GL}}
\def\Hom{{\rm Hom}}

\def\nt{\N^{\times}}
\def\N{{\mathbb N}}

\def\scal2{{\mathscr S}}
\def\HH{\text{HH}}
\def\strong{{\bf S\,}}

\begin{document}

\title{Hochschild homology, trace map  and $\zeta$-cycles}

\author{Alain Connes}
\address{Coll\`ege de France and IHES, France}
\email{alain@connes.org}

\author{Caterina Consani}
\address{Department of Mathematics, Johns Hopkins University,
Baltimore MD 21218, USA}
\email{cconsan1@jhu.edu}
\thanks{The second author is partially supported by the Simons Foundation collaboration grant n. 691493.}

\date{\today and, in revised form, \today.}

\dedicatory{}


\begin{abstract}In this paper we   consider two spectral realizations of the zeros of the Riemann zeta function. The first one involves all non-trivial (\ie non-real) zeros and  is expressed in terms of a Laplacian intimately related to the prolate wave operator. The  second spectral realization  affects only the critical zeros and it is  cast in terms of sheaf cohomology. The novelty is that the base space is the Scaling Site playing the role of the parameter space for the $\zeta$-cycles and encoding their stability by coverings. 
\end{abstract}

\maketitle

\section{Introduction}

In this paper we give a Hochschild  homological interpretation of  the  zeros of the Riemann zeta function.  The root of this result is in the recognition that the map 
\[
(\cE f)(u)=u^{1/2}\sum_{n>0} f(nu)
\]
which is defined on a suitable subspace  of the linear space of complex-valued even Schwartz functions  on the real line, is a trace  in Hochschild homology, if one brings in the construction the projection $\pi:\A_\Q\to \Q^\times\backslash \A_\Q$ from the rational ad\`eles to the ad\`ele classes (see Section \ref{sectgeom}). In this paper, we shall  consider two spectral realizations of the zeros of the Riemann zeta function. The first one involves all non-trivial (\ie non-real) zeros and  is expressed in terms of a Laplacian intimately related to the prolate wave operator (see Section \ref{sectlapla}). The  second spectral realization is sharper inasmuch as it affects only the critical zeros. The main players   are here the $\zeta$-cycles introduced in \cite{ccspectral}, and  the Scaling Site \cite{CCscal1} as their parameter space, which encodes their stability  by coverings. 
 The $\zeta$-cycles give the theoretical geometric explanation for the striking  coincidence between the low lying spectrum of a perturbed spectral triple therein introduced (see \cite{ccspectral}), and the low lying (critical) zeros
of the Riemann zeta function. The definition of a $\zeta$-cycle  derives, as a by-product,  from  scale-invariant Riemann sums  for complex-valued functions  on the real half-line $[0,\infty)$  with vanishing integral. For any $\mu\in\R_{>1}$, one implements the  linear (composite) map 
\[
\Sigma_\mu \cE: \sr0\to L^2(C_\mu)
\]
from  the Schwartz space $\sr0$ of real valued even functions  $f$ on the real line, with $f(0)=0$, and vanishing  integral, to the Hilbert space $L^2(C_\mu)$ of square integrable functions on the circle $C_\mu=\R_+^*/\mu^\Z$  of length $L=\log \mu$, where 
\[
(\Sigma_\mu g)(u):=\sum_{k\in\Z} g(\mu^ku).
\]
 The map $\Sigma_\mu$ commutes with the scaling action $\R^*_+\ni \lambda\mapsto f(\lambda^{-1}x)$  on functions, while $\cE$ is invariant under a normalized scaling action on $\sr0$.  In this set-up  one  has
 \begin{definition*} A {\bf $\zeta$-cycle} is  a circle $C$ of length $L=\log \mu$  whose Hilbert space $L^2(C)$ contains  $\Sigma_\mu \cE(\sr0)$ as a non dense subspace.
\end{definition*}
 
 Next result  is known (see \cite{ccspectral} Theorem 6.4)
\begin{theorem}\label{spectralreal} The following facts hold
\begin{enumerate}
    \item[(i)] The spectrum of the scaling action  of  $\R_+^*$ on the orthogonal space  to $\Sigma_\mu \cE(\sr0)$ in $L^2(C_\mu)$  is contained in the set of the imaginary parts of the zeros of the Riemann zeta function $\zeta(z)$ on the critical line $\Re(z)=\frac 12$.
\item[(ii)] Let $s>0$ be a real number such that $\zeta(\frac 12+is)=0$. Then any circle $C$ whose length is an integral multiple of $\frac{2\pi} s$ is a $\zeta$-cycle, and the spectrum of the action of $\R_+^*$ on  $(\Sigma_\mu \cE(\sr0))^\perp$ contains $s$.
\end{enumerate}
\end{theorem}

Theorem \ref{spectralreal} states that for a countable and dense set of values of $L\in\R_{>0}$, the Hilbert spaces $\cH(L):=(\Sigma_\mu \cE(\sr0))^\perp$ are non-trivial  and, more importantly, that as $L$ varies in that set,  the spectrum of the scaling action of $\R_+^*$ on the family of the $\cH(L)$'s is  the set $Z$ of imaginary parts of critical zeros of the Riemann zeta function. In fact, in view of the proven stability of $\zeta$-cycles under coverings, the same element of $Z$ occurs infinitely many times in the family of the $\cH(L)$'s. 
This stability under coverings displays the Scaling Site $\scal2=\rnt$ as the natural parameter space for the $\zeta$-cycles. In this paper, we show (see Section \ref{secttech}) that after organizing the family  $\cH(L)$ as a sheaf over $\scal2$  and using sheaf cohomology, one obtains a  spectral realization of  critical zeros of the Riemann zeta function. The key operation in the construction of the relevant arithmetic sheaf  is given by  the  action of the multiplicative monoid $\N^\times$  on the sheaf  of smooth sections of the bundle $L^2$ determined by the family of Hilbert spaces $L^2(C_\mu)$, $\mu={\exp L}$,  as $L$ varies in  $(0,\infty)$. For each $n\in\N^\times$ there is a canonical covering map $C_{\mu^n}\to C_\mu$,  where the action of $n$  corresponds to the operation of sum on the preimage of a point in $C_\mu$ under the covering.  This  action turns the  (sub)sheaf of smooth sections vanishing at $L=0$ into a sheaf $\cL^2$ over  $\scal2$. The family of  
 subspaces $\Sigma_\mu \cE(\sr0)\subset L^2(C_\mu)$  generates a closed subsheaf $\overline{\Sigma\cE}\subset \cL^2$    and one  then considers the cohomology of the related  quotient sheaf  $\cL^2/\overline{\Sigma\cE}$.  In view of the property of $\R_+^*$-equivariance under scaling, this construction  determines a spectral realization of  critical zeros of the Riemann zeta function, also taking care of eventual  multiplicities. Our main result is the following
 
  \begin{theorem}\label{spectralreal1intro} The  cohomology $H^0(\scal2,\cL^2/\overline{\Sigma\cE})$ endowed with the induced canonical action of $\R_+^*$ is isomorphic to the  spectral realization of  critical zeros of the Riemann zeta function, given by the action of $\R_+^*$,  via multiplication with $\lambda^{is}$, on the quotient of the Schwartz space $\cS(\R)$ by the closure of the ideal generated by multiples of $\zeta\left(\frac 12 +is\right)$. 	
 \end{theorem}

This paper is organized as follows. Section \ref{sectE} recalls  the main role played by the (image of the) map $\cE$
in the study of the spectral realization of the critical zeros of the Riemann zeta function.
In Section \ref{sectgeom} we show the identification of the Hochschild homology HH$_0$ of the 
noncommutative space $\Q^\times\backslash \A_\Q$ with the coinvariants for the action of $\Q^\times$ on the Schwartz algebra, using the (so-called) ``wrong way'' functoriality map $\pi_!$
associated to the projection $\pi:\A_\Q\to \Q^\times\backslash \A_\Q$. We also stress the relevant fact that the Fourier transform on ad\`eles becomes canonical after passing to $\HH_0$ of the ad\`ele class space of the rationals. The key Proposition \ref{propmapE1} describes the invariant part of  such $\HH_0$ as the space of even Schwartz functions on the real line and identifies the trace map with the map $\cE$. 
Section \ref{sectlapla} takes care of the two vanishing conditions implemented in the definition of  $\cE$ and introduces the  operator $\Delta=H(1+H)$ ($H$ being the generator of the scaling action of $\R^*_+$ on $\cS(\R)^{\text{ev}}$) playing the role of the Laplacian and intimately related to the prolate operator. 
Finally, Section \ref{secttech} is the main technical section of this paper since it contains the proof of Theorem \ref{spectralreal1intro}. 

 \section{The map $\cE$ and the zeros of the zeta function}\label{sectE}
  The ad\`ele class space of the rationals $\Q^\times\backslash \A_\Q$  is the natural geometric framework to understand the Riemann-Weil explicit formulas for $L$-functions as a trace formula \cite{Co-zeta}. The essence of this result lies mainly in the delicate computation of the  principal values involved in the distributions appearing in the  geometric (right-hand) side of the semi-local trace formula of \opcit (see Theorem 4 for the notations)
  \[
  \text{Trace}(R_\Lambda U(h))=h(1)\int_{\lambda\in C_S\atop |\lambda|\in[\Lambda^{-1},\Lambda]}d^*\lambda+\sum_{\nu\in S}\int'_{\Q_\nu^*}\frac{h(u^{-1})}{|1-u|}d^*u+o(1)\qquad \text{for $\Lambda\to\infty$}
  \]
  (later recast in the softer context of \cite{Meyer}). There is a rather simple analogy related to the spectral (left-hand) side of the explicit  formulas for a global field $K$ (see \cite{CRH} Section~2 for the notations)
  \[
  \hat h(0)+\hat h(1) - \sum_{\chi\in \widehat{C_{K,1}}}\sum_{\rho\in Z_{\tilde\chi}}\hat h(\tilde\chi,\rho) = \sum_\nu\int'_{K_v^*}\frac{h(u^{-1})}{|1-u|}d^*u
  \]
  which may help one to realize how the sum over the zeros of the zeta function appears. Here this relation is simply explained.
  Given a complex valued polynomial $P(x)\in \C[x]$, one may identify the set of its zeros  as the spectrum of the endomorphism $T$ of multiplication by the variable $x$ computed in the quotient algebra $\C[x]/(P(x))$. It is well known that the matrix of $T$, in the basis of powers of $x$, is the companion matrix of $P(x)$. Furthermore, the trace of its powers, readily computed from the diagonal terms of powers of the companion matrix in terms of the coefficients of $P(x)$, gives the Newton-Girard formulae.\footnote{This is an efficient way to find the power sum of roots of $P(x)$ without actually finding the roots explicitly.  Newton's identities supply the calculation via a recurrence relation with known coefficients.}
 If one transposes this result to the case of the Riemann zeta function $\zeta(s)$, one sees that the multiplication by $P(x)$ is replaced here with the  map 
 \begin{equation}\label{mapE}
 \cE(f)(u):=u^{1/2}\sum_{n=1}^\infty f(nu),
 \end{equation} 
 while the role of $T$ (the multiplication by the variable) is played by the scaling operator $u\partial_u$. These statements may become  more evident if one brings in the Fourier transform.
 Indeed, let $f\in \cS(\R)^{\rm ev}$ be an even Schwartz function  and let $w(f)(u)=u^{1/2}f(u)$ be the unitary identification  of $f$ with a function in $L^2(\R_+^*,d^*u)$, where  $d^*u:=du/u$ denotes the Haar measure. Then, by composing $w$ with the  (multiplicative) Fourier transform $\fourier: L^2(\R_+^*,d^*u)\to L^2(\R)$,  $\fourier(h)(s)=\int_{\R_+^*}h(u)u^{-is}d^*u$ one obtains 
 \[
 \fourier(w(f))=\psi, \qquad   
\psi(z)=\int_{\R_+^*}f(u)u^{\frac 12-iz}d^*u.
\]
The function $\psi(z)$ is holomorphic in the complex half-plane $\mathfrak H=\{z\in\C~|~\Im(z)>-\frac 12\}$ since $f(u)=O(u^{-N})$ for $u\to \infty$. Moreover, for $n\in\N$, one has
\[
\int_{\R_+^*}u^{1/2}f(nu)u^{-iz}d^*u=n^{-1/2+iz}\int_{\R_+^*}v^{1/2}f(v)v^{-iz}d^*v.
\]
In the region $\Im(z)>\frac 12$ one derives, by applying Fubini theorem, the following equality
\[
\int_{\R_+^*}{\sum_{n=1}^\infty} u^{1/2}f(nu)u^{-iz}d^*u= \left({\sum_{n=1}^\infty} n^{-1/2+iz}\right)\int_{\R_+^*}v^{1/2}f(v)v^{-iz}d^*v.
\]
Thus, for all $z\in\C$ with $\Im(z)>\frac 12$ one obtains 
\begin{equation}\label{mapEzeta}
	\int_{\R_+^*}\cE(f)(u)u^{-iz}d^*u=\zeta(\frac 12-iz)\psi(z).
\end{equation}
If one assumes now that the Schwartz function $f$ fulfills  $\int_\R f(x)dx=0$, then $\psi(\frac i2)=0$. Both  sides of \eqref{mapEzeta} are holomorphic functions in  $\mathfrak H$: for the integral on the left-hand side, this can be seen by using the estimate $\cE(f)(u)=O(u^{1/2})$ that follows from the Poisson formula. This proves that \eqref{mapEzeta} continues to hold also in the complex half-plane $\mathfrak H$.  Thus one sees that  the zeros of $\zeta(\frac 12-iz)$  in the strip $\vert \Im(z)\vert <\frac 12$ are the common zeros of all functions
 \begin{equation*}\label{mapErange}
 \fourier(\cE(f))(z), \qquad f \in \cS(\R)^{\rm ev}_1:=\{ f \in \cS(\R)^{\rm ev}~|~  \int_\R f(x)dx=0\}.
 \end{equation*} 
One may eventually select the even Schwartz function $f(x)=e^{-\pi  x^2} \left(2 \pi  x^2-1\right)$ to produce a specific instance where the zeros of $\fourier(\cE(f))$ are exactly the non-trivial zeros of $\zeta(\frac 12-iz)$, since in this case  $\psi(z)=\frac{1}{4} \pi ^{-\frac{1}{4}+\frac{i z}{2}} (-1-2 i z) \Gamma \left(\frac{1}{4}-\frac{i z}{2}\right)$.

\section{Geometric interpretation}\label{sectgeom}
In this section we  continue the study of the map $\cE$ with the goal to achieve a geometric understanding of it. This is obtained by bringing in the construction the ad\`ele class space of the rationals, whose role  is that to grant  for the replacement, in \eqref{mapE}, of  the summation over the monoid $\nt$  with the summation over the group $\Q^\times$. Then, up to the factor $u^{1/2}$, $\cE$ is understood as  the composite $\iota^*\circ \pi_!$, where the map $\iota:\Q^\times\backslash\A^*_\Q/\hatz\to \Q^\times\backslash \A_\Q/\hatz$ is the inclusion of id\`ele classes in ad\`ele classes and  $\pi:\A_\Q/\hatz\to \Q^\times\backslash \A_\Q/\hatz$ is induced by the projection $\A_\Q\to \Q^\times\backslash \A_\Q$. We shall discuss the following diagram  
\begin{equation*}\label{mapE1}
\xymatrix{ Y=\A_\Q/\hatz  \ar[d]^\pi &  \\
X=\Q^\times\backslash \A_\Q/\hatz  & \Q^\times\backslash \A^\times_\Q/\hatz=\R_+^*\ar[l]_{\iota} }
           \end{equation*}
 The conceptual understanding of the map $\pi_!$ uses Hochschild homology of noncommutative algebras.  
 
 We recall that the space of ad\`ele classes \ie the quotient  $\Q^\times\backslash \A_\Q$ is encoded algebraically by the cross-product algebra
\begin{equation*}\label{crossp}
 \cA:=\cS(\A_\Q)\rtimes  \Q^\times.
 \end{equation*}
 The Schwartz space $\cS(\A_\Q)$ is acted upon by (automorphisms of) $\Q^\times$ corresponding to the scaling action of  $\Q^\times$ on rational ad\`eles. An element of $\cA$ is written symbolically as a finite sum
 \[
 \sum a(q)U(q), \qquad a(q)\in \cS(\A_\Q).
 \]
 From the inclusion of algebras $\cS(\A_\Q)\subset \cS(\A_\Q)\rtimes  \Q^\times= \cA$ one derives  a corresponding morphism of Hochschild homologies 
 \begin{equation*}\label{crossp1}
\pi_!: \text{HH}(\cS(\A_\Q))\longrightarrow \text{HH}(\cA).
\end{equation*}
 Here, we use the shorthand notation $\text{HH}(A):=\text{HH}(A,A)$ for the Hochschild homology of an algebra $A$ with coefficients in the bimodule $A$. In  noncommutative geometry, the vector space of differential forms of degree $k$ is replaced by the Hochschild homology $\HH_k(A)$. If the algebra $A$ is commutative and for $k=0$,  $\HH_0(A)=A$, so that  $0$-forms are identified with functions. Indeed, the Hochschild boundary map
 \[
 b:A^{\otimes\, 2}\to A, \qquad b(x\otimes y)=xy-yx
 \]
 is identically zero when the algebra $A$ is commutative. This result does not hold when $A=\cA$, since   $\cA=\cS(\A_\Q)\rtimes  \Q^\times$ is  no longer commutative. It is therefore meaningful to bring in the following
 \begin{proposition}\label{h0coinv} The kernel of  $\pi_!: \HH_0(\cS(\A_\Q))\to \HH_0(\cA)$ is the $\C$-linear span $E$ of functions  $f-f_q$, with $f\in \cS(\A_\Q)$, $q\in \Q^\times$, and where we set $f_q(x):=f(qx)$. 	
 \end{proposition}
 \proof For any $f,g\in \cS(\A_\Q)$ and  $q\in \Q^\times$ one has 
 \begin{equation}\label{crossp2}
 fg-(fg)_q=fg-U(q)fg U(q^{-1})=xy-yx, \qquad x:=f U(q^{-1}), \ y:=U(q)g.
 \end{equation}
One knows (\cite{Miyazaki} Lemma 1) that any function $f\in \cS(\R)$ is a product of two elements of $\cS(\R)$. Moreover, an element of the Bruhat-Schwartz space $\cS(\A_\Q)$ is  a finite linear combination of functions of the form $e\otimes f$, with $e^2=e$. Thus any $f\in \cS(\A_\Q)$ can be written as a finite sum of products of two elements of  $\cS(\A_\Q)$, so that \eqref{crossp2} entails $f-f_q\in \ker \pi_!$. Conversely, let  $f\in \ker \pi_!$. Then there exists a finite number of pairs $x_i, y_i \in \cA$ such that $f=\sum [x_i,y_i]$. Let $P: \cA\to \cS(\A_\Q)$ be the projection on the coefficient $a(1)$ of $U(1)=1$ \ie
\[
P\left( \sum a(q)U(q)\right) :=a(1).
\]
Then  $f=\sum [x_i,y_i]$ implies $f=\sum P( [x_i,y_i])$. We shall prove that for any pair $x, y \in \cA$ one has $P( [x,y])\in E$. Indeed, one has 
\[
x=\sum a(q)U(q), \ y=\sum b(q') U(q'),\ [x,y]=\sum (a(q)b(q')_q-b(q')a(q)_{q'})U(qq')
\]
so that 
\[
P( [x,y])=\sum (a(q)b(q^{-1})_q-b(q^{-1})a(q)_{q^{-1}}).
\]
This projection  belongs to $E$ in view of the fact that
\[
a(q)b(q^{-1})_q-b(q^{-1})a(q)_{q^{-1}}=h_q-h, \qquad h=b(q^{-1})a(q)_{q^{-1}}.
\]
This completes the proof.\endproof 

Proposition \ref{h0coinv} shows that the image of  $\pi_!: \HH_0(\cS(\A_\Q))\to \HH_0(\cA)$ is the space of coinvariants for the action of $\Q^\times$ on $\cS(\A_\Q)$, \ie the quotient of $\cS(\A_\Q)$ by the subspace $E$.\newline
 An important point now to remember is  that the Fourier transform becomes canonically defined on the above quotient. Indeed, the definition of the Fourier transform on ad\`eles depends on the choice of a non-trivial character $\alpha$ on the additive, locally compact  group $\A_\Q$, which is trivial on the subgroup $\Q\subset \A_\Q$. It is  defined as follows
 \begin{equation*}\label{fourier1}
\fourier_\alpha(f)(y):= \int_\R f(x)\alpha(xy) dx.
\end{equation*}
 The space of characters of the compact group $G=\A_\Q/\Q$ is one dimensional as a $\Q$-vector space, thus any non-trivial character $\alpha$ as above is of the form $\beta(x)=\alpha(qx)$, so that 
 \begin{equation*}\label{fourier2}
\fourier_\beta(f)(y):= \int_\R f(x)\alpha(qxy) dx=\fourier_\alpha(f)_q(y).
\end{equation*}
Therefore, the difference $\fourier_\beta-\fourier_\alpha$ vanishes on the quotient of $\cS(\A_\Q)$ by $E$ and this latter space is preserved by $\fourier_\alpha$ since  $\fourier_\alpha(f_q)=\fourier_\alpha(f)_{q^{-1}}$. 

\subsection{$\HH$, Morita invariance and the trace map} Let us recall that given an algebra $A$, the trace map 
\[
\tr : M_n(A)\to A, \qquad \tr((a_{ij}):=\sum a_{jj}
\]
induces an isomorphism in degree zero Hochschild homology 
which extends to higher degrees. If $A$ is a convolution algebra of the \'etale groupoid of an equivalence relation $\cR$ with countable orbits on a space $Y$, and  $\pi:Y\to Y/\cR$ is the quotient map, the trace map takes the following form 
 \begin{equation}\label{morita0}
\tr(f)(x):=\sum_{\pi(j)=x} f(j,j).
\end{equation}
The trace induces a map on $\HH_0$ of the function algebras, provided one takes care of the convergence issue when the size of equivalence classes is infinite. If the relation $\cR$  is associated with the orbits of the free  action of a discrete  group $\Gamma$ on  a locally compact space $Y$, the convolution algebra is the cross product of the algebra of functions on $Y$ by the discrete  group $\Gamma$. In this case, the \'etale groupoid is $Y\rtimes \Gamma$, where the source and range maps are given {\it resp.} by $s(y,g)=y$ and $r(y,g)=gy$. The elements of the convolution algebra are functions $f(y,g)$ on $Y\rtimes \Gamma$. The diagonal terms  in \eqref{morita0} correspond to the elements of $Y\rtimes \Gamma$ such that $s(y,g)=r(y,g)$, meaning that $g=1$ is the neutral element of $\Gamma$, since the action of $\Gamma$ is assumed to be free. Then, the trace map  is
 \begin{equation*}\label{morita}
\tr((f)(x)=\sum_{\pi(y)=x} f(y,1).
\end{equation*}
This sum is meaningful on the space of the proper orbits of $\Gamma$. For a lift $\rho(x)\in Y$, with $\pi(\rho(x))=x$ the trace reads as \begin{equation}\label{morita1}
\tr(f)(x)=\sum_{g\in \Gamma} f(g\rho(x),1).
\end{equation}
In the case of    $Y=\A_\Q$   acted upon by $\Gamma=\Q^\times$, the  proper orbits are parameterized by the id\`ele classes and this space embeds in the ad\`ele classes by means of the inclusion \[
\iota : \Q^\times\backslash \A^\times_\Q\to \Q^\times\backslash \A_\Q.
\]
We identify the id\`ele class group $C_\Q=\Q^\times\backslash \A^\times_\Q$  with  $\hatz\times \R_+^*$, using the canonical exact sequence affected by the modulus
\[
1\to \hatz\to C_\Q \stackrel{\Mod}{\longrightarrow} \R_+^*\to 1.
\]
There is a natural section $\rho:C_\Q\to \A^\times_\Q$ of the quotient map,  given by the canonical inclusion 
\[
\hatz\times \R_+^*\subset \A_\Q^f\times \R =\A_\Q.
\]
Next, we focus on the $\hatz$-invariant part  of $\cS(\A_\Q)$. Then, with the notations of Proposition~\ref{h0coinv} we have

\begin{lemma}\label{propmapE}  
The following facts hold\begin{enumerate}
    \item[(i)] Let $h\in\cS(\A_\Q)^{\hatz}$, then there exists $ f\in \cS(\R)^{\rm ev}$ with $h-1_{\hat \Z}\otimes f\in E^{\hatz}$.
\item[(ii)] Let $ f\in \cS(\R)^{\rm ev}$ and $\tilde f=(1_{\hat \Z}\otimes f)U(1)\in S(\A_\Q)^{\hatz}\rtimes \Q^\times$, then one has 
\begin{equation}\label{tracemap}
\tr(\tilde f)(u)=2\sum_{n\in \nt} f(nu)\qquad \forall u\in \R_+^*.
\end{equation}
\item[(iii)] Let $ f\in \cS(\R)^{\rm ev}$, then $1_{\hat \Z}\otimes f\in E^{\hatz}\iff f=0$.
\end{enumerate}
 \end{lemma} 
\proof (i)~By definition, the elements of the Bruhat-Schwartz space $\sarr$  are finite linear combinations of functions on $\A_\Q$ of the form  ($S\ni \infty$ is a finite set of places)
\[
f=\otimes f_v, \qquad f_v=1_{\Z_v}~  \forall v \notin S, \quad  f_{\infty}\in \cS(\R), \quad f_p\in \cS(\Q_p)\quad \forall p\in S\setminus \infty,
\]
where $\cS(\Q_p)$ denotes the space of locally constant functions with compact support. An element of 
$\cS(\Q_p)$ which is $\Z_p^*$-invariant is a finite linear combination of  characteristic functions $(1_{\Z_p})_{p^n}(x):=1_{\Z_p}(p^n x)$. Thus an element  $h\in \cS(\A_\Q)^{\hatz}$ is   a finite linear combination of functions of the form
\[
f=\otimes f_v, \qquad f_v=1_{\Z_v}\quad \forall v \notin S, \quad  f_{\infty}\in \cS(\R), \quad f_p=(1_{\Z_p})_{p^{n_p}}\quad \forall p\in S\setminus \infty
\]
With $q=\prod p^{-n_p}$ one has with $\ell(x):=f(qx)$,  $\ell=1_{\hat \Z}\otimes g$, $\ell -f\in E^{\hatz}$ and the replacement of $g$ with its even part $\frac 12 (g(x)+g(-x))$ does not change the class of $f$ modulo $E^{\hatz}$.\newline
(ii)~Let $f\in \cS(\R)^{\rm ev}$ and $\tilde f=(1_{\hat \Z}\otimes f)U(1)\in S(\A_\Q)^{\hatz}\rtimes \Q^\times$. For $u\in \R_+^*$ the lift $\rho(u)\in \cS(\A_\Q)$ is  $(1, u)\in \hatz\times \R_+^*\subset \A_\Q^f\times \R =\A_\Q$, and by applying \eqref{morita1} one has
\begin{equation}\label{trmap}
\tr(\tilde f)(u)=\sum_{q\in \Q^\times} \tilde f(q \rho(x),1)=
\sum_{q\in \Q^\times} (1_{\hat \Z}\otimes f)(q,qu)=2\sum_{n\in \nt} f(nu)\quad \forall u\in \R_+^*.
\end{equation}
(iii)~Let $ f\in \cS(\R)^{\rm ev}$ such that $1_{\hat \Z}\otimes f\in E^{\hatz}$. Let $\tilde f=(1_{\hat \Z}\otimes f)U(1)\in S(\A_\Q)^{\hatz}\rtimes \Q^\times$. By Proposition \ref{h0coinv} the Hochschild class in $\HH_0(\cA)$ of  $\tilde f$ is zero, thus $\tr(\tilde f)=0$. It follows from \eqref{tracemap}  
that $\cE(f)(u)=0$ $\forall u\in \R_+^*$. Then \eqref{mapEzeta} implies that the function $\psi(z)=\int_{\R_+^*}f(u)u^{\frac 12-iz}d^*u$ is well defined in the half-plane $\Im(z)  >\frac 12$ where it vanishes identically, thus $f=0$. The converse of the statement is obvious.
\endproof

The next statement complements Proposition~\ref{h0coinv}, with a description of the  range of $\pi_!: \HH_0(\cS(\A_\Q)^{\hatz})\to \HH_0(\cA)^{\hatz}$;  it also shows that the map $\cE(f)(u)=u^{1/2}\sum_{n=1}^\infty f(nu)$ coincides, up to the factor $\frac{u^{1/2}}{2}$, with the trace map \eqref{trmap}. We keep the notations of Lemma~\ref{propmapE}

\begin{proposition}\label{propmapE1}
The map $\cS(\R)^{\rm ev}\to \pi_!\left(\HH_0(\cS(\A_\Q)^{\hatz})\right)$, $f\mapsto \tilde f$ is an isomorphism, this means that  $\pi_!\left(\HH_0(\cS(\A_\Q)^{\hatz})\right)$ is determined by the images of the elements of the subalgebra $1_{\hat \Z}\otimes\cS(\R)^{\rm ev}\subset \cS(\A_\Q)^{\hatz}$. Furthermore, one has the identity
\begin{equation*}\label{hhmapE}
\cE(f)(u)=\frac{u^{1/2}}{2}\ \tr(\tilde f)(u)\quad f\in \cS(\R)^{\rm ev},~\forall u\in \R_+^*.
\end{equation*}
\end{proposition}
\proof The first statement follows from Lemma \ref{propmapE} (i) and (iii). The second statement from (ii) of the same lemma. \endproof

  \section{The  Laplacian $\Delta=H(1+H)$}\label{sectlapla}
  
  This section describes the spectral interpretation of the squares of non-trivial zeros of the Riemann zeta function in terms of a suitable Laplacian. It also shows the relation between this Laplacian and the prolate wave operator.
  
  \subsection{The vanishing conditions} 
  One starts with the exact sequence
\begin{equation*}\label{fonction1}
0\to \cS(\A_\Q)_1 \to \cS(\A_\Q)\stackrel{\epsilon}{\to} \C(1)\to 0
\end{equation*}
associated to the kernel of the $\Q^\times$-invariant linear functional $\epsilon(f)=\int_{\A_\Q} f(x)dx\in  \C(1)$. By implementing in the above sequence the evaluation $\delta_0(f):=f(0)$, one obtains the exact sequence
\begin{equation*}\label{fonction0}
0\to \cS(\A_\Q)_0 \to \cS(\A_\Q)\stackrel{(\delta_0,\epsilon)}{\longrightarrow} \C(0)\oplus \C(1)\to 0.
\end{equation*}
The next lemma shows that both $\cS(\A_\Q)_0$ and $\cS(\A_\Q)_1$ have a description in terms of the ranges of two related differential operators. For simplicity of exposition, we restrict our discussion to the $\hatz$-invariant parts of these function spaces.

\begin{lemma} \label{quotbruhatschwartz} Let $H:\cS(\A_\Q)\to \cS(\A_\Q)$, $H:=x\partial_x$ be the generator of the scaling action of $1\times \R_+^*\subset \GL_1(\A_\Q)$. Then one has
\begin{enumerate}
\item[(i)] $H$  commutes with the action of $\GL_1(\A_\Q)$ and  restricts 
to $\cS(\A_\Q)^{\hatz}$.
\item[(ii)] $(1+H)$  induces an isomorphism $\cS(\R)^{\rm ev}\to \cS(\R)_1^{\rm ev}$.
\item[(iii)] $H(1+H)$  induces an isomorphism $\cS(\R)^{\rm ev}\to \cS(\R)_0^{\rm ev}$.
\end{enumerate}
\end{lemma}
\proof  (i)~follows since $\GL_1(\A_\Q)$ is abelian, thus $H$  commutes with the action of $\GL_1(\A_\Q)$.\newline
(ii)~ The functional $f\mapsto \epsilon(f)=\int_{\A_\Q}  f(x)dx$ vanishes on the range of $1+H$ since
\[
\epsilon((1+H)f)=\int_{\A_\Q} \partial(xf)dx=0,
\]
 thus the range of $(1+H)$ is contained in $\cS(\R)^{\rm ev}_1$. The equation  $Hf+f=0$ implies that $xf(x)$ is constant, hence  $f=0$, for $f\in \cS(\R)$. Thus $(1+H):\cS(\R)^{\rm ev}\to \cS(\R)_1^{\rm ev}$ is injective. Let now  $f\in\cS(\R)^{\rm ev}_1$ and let $\widehat f$ be its Fourier transform. Then $\widehat f \in \cS(\R)^{\rm ev}$ and $\widehat f(0)=0$. The function $g(x):=\widehat f(x)/x$, $g(0):=0$ is smooth, and $g\in\cS(\R)^{\rm odd}$, while the function 
 \[
 h(x):=\int_{-\infty}^x g(y)dy
 \]
fulfills $h\in\cS(\R)^{\rm ev}$ and $\partial_x h=g$. One has $Hh=\widehat f$, thus $(-1-H)\widehat h= f$. This shows that $(1+H):\cS(\R)^{\rm ev}\to \cS(\R)_1^{\rm ev}$ is also surjective. \newline
(iii)~Since the evaluation $f\mapsto f(0)$ is invariant under scaling, it vanishes on the range of $H$. Similarly, one sees that the functional $f\mapsto \int_{\A_\Q}  f(x)dx$ vanishes on the range of $1+H$. Thus the range of $H(1+H)$ is contained in $\cS(\R)_0$. The equation $Hf=0$ implies that the function $f$ is constant and hence $f=0$, for $f\in \cS(\R)$. Similarly $Hf+f=0$ implies that $xf(x)$ is constant and hence $f=0$ for $f\in \cS(\R)$. Thus $H(1+H): \cS(\R)\to \cS(\R)_0$ is injective. Let now $f\in\cS(\R)^{\rm ev}$  with $f(0)=0$. Then the function $g(x):=f(x)/x$, $g(0):=0$, is smooth, $g\in\cS(\R)^{\rm odd}$ and there exists a unique $h\in\cS(\R)^{\rm ev}$ such that $\partial_x h=g$. One has $Hh=f$ so that  $(-1-H)\widehat h=\widehat f$. Thus if $\widehat f(0)=0$ one has $\widehat h(0)=0$ and there exists $k\in\cS(\R)^{\rm ev}$ with $Hk=\widehat h$. Then $-(1+H)\widehat k=h$ and $H(1+H)\widehat k=-f$. This shows that $H(1+H):\cS(\R)^{\rm ev}\to \cS(\R)_0^{\rm ev}$ is surjective and an isomorphism.\endproof 

\subsection{The Laplacian $\Delta=H(1+H)$ and its spectrum}

This section is based on the following heuristic dictionary suggesting a parallel  between some classical notions in Hodge theory on the left-hand side,  and their counterparts in noncommutative geometry, for the ad\`ele class space of the rationals. The notations  are inclusive of those of Section \ref{sectgeom}
\begin{large}
\begin{center}
\begin{tabular}{l | l}
& \\
Algebra of functions & Cross-product by $\Q^\times$ \\
Differential forms  & Hochschild homology\\
Star operator $\star$ & $\iota\ \times $ $\fourier$ \\ 
 Differential $d$  & Operator $H$  \\
$\delta:= \star d \star$ & Operator $1+H$   \\
Laplacian  & $\Delta:=H(1+H)$ \\
&
\end{tabular}
\end{center}
\end{large}
Next Proposition is a variant of the spectral realization in  \cite{CCM2,CMbook}.

  \begin{proposition} \label{laplacian} The following facts hold
  \begin{enumerate}
      \item[(i)] The trace map $\tr$ commutes with $\Delta=H(1+H)$ and the range of $\tr\circ \Delta$ is contained in the strong Schwartz space $\strong(\R_+^*):=\,\cap_{\beta \in \R}\,\mu^\beta \cS(\R_+^*)$, with
      $\mu$ denoting the Modulus.
  \item[(ii)] The spectrum of $\Delta$ on the quotient of  $\strong(\R_+^*)$ by the closure of the range of  $\tr\circ \Delta$ is the set (counted with possible multiplicities)
  \[
  \left\{\Big(z-\frac 12\Big)^2-\frac 14\mid  z\in \C\setminus \R, ~ \zeta(z)=0\right\}.
  \]
  \end{enumerate}
  \end{proposition}
\proof (i)~The trace map of \eqref{trmap}  commutes with $\Delta$. By Lemma \ref{quotbruhatschwartz} (iii) the range of $\Delta$ is $\cS(\R)_0^{\rm ev}$ thus the range of $\cE\circ(H(1+H))$ is contained in  $\strong(\R_+^*)$ (see \cite{CMbook}, Lemma 2.51).\newline
(ii)~By construction,  $\strong(\R_+^*)$ is  the intersection, indexed by compact intervals $J\subset \R$, of  
the spaces $\cap_{\beta \in J}\,\mu^\beta \cS(\R_+^*) $. The Fourier transform  
\begin{equation*}\label{fouriermellin}
\fourier(f)(z):=\int_{\R_+^*}f(u)u^z d^*u	
\end{equation*}
defines an isomorphism of these function spaces with the Schwartz spaces $\cS(I)$, labeled by vertical strips $I:=\{z\in \C\mid \Re(z)\in J\}$,
of holomorphic functions $f$ in $I$ with $p_{k,m}(f)<\infty$ for all $k,m\in\N$ where 
\begin{equation*}\label{seminormspkm}
p_{k,m}(f)=\sum_{n=0}^m\,\frac{1}{n!}\,\sup_I (1+|z|)^k\cdot|\partial^n
f(z)|\,.
\end{equation*}
These norms define the Fr\'echet topology on $\cS(I)$. \newline
It then follows from Lemma 4.125 of \cite{CMbook} that, for $I$ sufficiently large, the quotient of $\cS(I)$  by the closure of the range of  $\tr\circ \Delta$ decomposes into a direct sum of finite-dimensional spaces associated to the projections $\Pi(N)$, $N\in \Z$ which fulfill the following properties
\begin{enumerate}
  \item Each $\Pi(N)$  is an idempotent.
  \item The sequence of $\Pi(N)$'s  is of tempered growth. \index{tempered growth}
   \item The rank of $\Pi(N)$  is $O(\log |N|)$ for $|N|\to \infty$.
   \item For any $f\in \cS(I)$ the sequence $f\,\Pi(N)$ is of rapid
   decay and \index{rapid decay}
   \begin{equation*}\label{opPiN6}
\sum_{N\in\Z}\,f\,\Pi(N)=f \qquad \forall f\in \cS(I)\,.
\end{equation*}
\end{enumerate}
This direct sum decomposition commutes with  $\Delta$ since both $\Pi(N)$ and the conjugate of $\Delta$ by the Fourier transform $\fourier$ are given by multiplication operators.  The conjugate of $H$  by  $\fourier$ is the multiplication by $-z$, so that the conjugate of $\Delta$ is the multiplication by $-z(1-z)$. 
The spectrum of $\Delta$ is the union of the spectra of the finite-dimensional operators $\Delta_N:=\Pi(N)\Delta=\Delta\Pi(N)$. By \cite{CMbook}, Corollary 4.118, and the proof of Theorem 4.116, the finite-dimensional range of $\Pi(N)$ is described by the evaluation of $f\in \cS(I)$  on the zeros $\rho\in Z(N)$ of the Riemann zeta function which are inside the contour $\gamma_N$, \ie by the map 
\begin{equation*}\label{evalmaprho}
\cS(I)\ni f\mapsto f|_{Z(N)}=(f^{(n_\rho)}(\rho))_{\rho\in Z(N)}\in \C^{(n_\rho)} 
\end{equation*}
where $\C^{(n_\rho)}$ denotes the space of dimension
$n_\rho$ of jets of order equal to the order $n_\rho$ of the zero $\rho$ of the zeta function. Moreover, the action of $\Delta_N$ is given by the matrix associated with the multiplication of $f\in \cS(I)$ by $-z(1-z)$:  this gives a triangular matrix whose diagonal is given by $n_\rho$ terms all equal to $-\rho(1-\rho)$. Thus 
  the spectrum of $\Delta$ on the quotient of  $\strong(\R_+^*)$ by the closure of the range of  $\tr\circ \Delta$ is the set (counted with multiplicities)
  \[
  \left\{\Big(\rho-\frac 12\Big)^2-\frac 14\mid \rho\in \C\setminus \R, \ \zeta(\rho)=0\right\}.
  \]
  \endproof 
  
\begin{corollary}\label{rhequiv} The spectrum of $\Delta$ on the quotient of the strong Schwartz space $\strong(\R_+^*)$ by the closure of the range of  $\tr\circ \Delta$ is negative if and only if the Riemann Hypothesis holds.	
\end{corollary}
\proof This follows from Proposition \ref{laplacian} and the fact that for $\rho \in \C$
\[
\Big(\rho-\frac 12\Big)^2 -\frac 14\leq 0 \iff \rho\in [0,1]\cup \frac 12 +i \R.
\]
 \endproof 
 
 \begin{remark}\label{prolate} The main interest of the above reformulation of the spectral realization of \cite{CCM2,CMbook} in terms of the Laplacian $\Delta$ is that the latter is intimately related to the prolate wave operator $W_\lambda$ that is shown in \cite{CMo} to be self-adjoint and have, for $\lambda=\sqrt 2$ the same UV spectrum as the Riemann zeta function. The relation between  $\Delta$ and 	$W_\lambda$ is that the latter is a perturbation of $\Delta$ by a multiple of the Harmonic oscillator.
 \end{remark}

 \section{Sheaves on the Scaling Site and  $H^0(\scal2,\cL^2/\overline{\Sigma\cE})$}\label{secttech}

Let $\mu\in\R_{>1}$ and  $\Sigma_\mu$    be  the linear map on functions $g: \R^*_+ \to \C$   of sufficiently rapid decay at $0$ and $\infty$  defined by 
\begin{equation}\label{sigmamu}
(\Sigma_\mu g)(u)=\sum_{k\in\Z}g(\mu^ku).
\end{equation}
We shall  denote with $\sr0$ the linear space of real-valued, even Schwartz functions $f \in \cS(\R)$ fulfilling the two conditions
$f(0) = 0 = \int_\R f(x)dx$. The map
\begin{equation}\label{Epsilon}
\cE: \sr0 \to \R,\qquad (\cE f)(u)= u^{1/2}\sum_{n>0}f(nu)
\end{equation}
is proportional to a Riemann sum for the integral of $f$.
 The following  lemma on scale invariant Riemann sums  justifies the pointwise ``well-behavior" of \eqref{Epsilon} (see \cite{ccspectral} Lemma~6.1)
 
 \begin{lemma}\label{fouriertruncated1}\; Let $f$ be a complex-valued function of bounded variation on $(0,\infty)$. Assume that $f$ is of rapid decay for $u\to \infty$,  $O(u^2)$ when $ u\to 0$, and  that $\int_0^{\infty}f(t)dt=0$.  Then the following properties hold 
 \begin{enumerate}
\item[(i)]  The function $(\cE f)(u)$ in \eqref{Epsilon} is well-defined pointwise, is 
$O(u^{1/2})$ when $u\to 0$, and of rapid decay for $u\to\infty$.
\item[(ii)] Let  $g=\cE(f)$, then the series \eqref{sigmamu}  is geometrically convergent, and defines a bounded and measurable  function on  $\R_+^*/\mu^{\Z}$.
\end{enumerate}
\end{lemma}

 We recall that a sheaf over  the Scaling Site $\scal2=\rnt$ is a sheaf of sets on $[0,\infty)$ (endowed with the euclidean topology)  which is equivariant for the action of the multiplicative monoid $\nt$ \cite{CCscal1}.  Since we work in characteristic zero we select as structure sheaf  of $\scal2$ the $\nt$-equivariant sheaf $\cO$ whose sections on an open set $U\subset [0,\infty)$ define the space of smooth, complex-valued  functions on $U$.
 The next proposition introduces two relevant sheaves of $\cO$-modules.
 
 \begin{proposition}\label{sheaves}
Let $L\in (0,\infty)$, $\mu=\exp L$, and $C_\mu=\R^*_+/\mu^\Z$. The following facts hold
 \begin{enumerate}
 \item[(i)] As $L$ varies in $(0,\infty)$, the pointwise multiplicative Fourier transform
 \begin{equation}\label{pointwisefourier}
	 \fourier:L^2(C_\mu)\to \ell^2(\Z), \qquad \fourier(\xi)(n)=L^{-\frac 12}\int_{C_\mu} \xi(u)u^{-\frac{2\pi i n}{L}}d^*u 
\end{equation}
defines an isomorphism between the family of Hilbert spaces $L^2(C_\mu)$  and the restriction to $(0,\infty)$ of the  trivial vector bundle $L^2=[0,\infty)\times \ell^2(\Z)$.
 \item[(ii)] The smooth sections  vanishing at $L=0$ of the vector bundle $L^2$   together with the linear maps 
 \begin{equation}\label{sigman}
 \sigma_n:L^2(C_{\mu^n})\to L^2(C_\mu), \qquad \sigma_n(\xi)(u)=\sum_{j=1}^n \xi(\mu^j u)
 \end{equation}
 define a  sheaf  $\cL^2$ of $\cO$-modules on $\scal2$.
 \item[(iii)] For $f\in \sr0$, the maps $\Sigma\cE(f)(L):= \Sigma_{\mu}\cE(f)$  as in \eqref{sigmamu} define  smooth (global) sections $\Sigma\cE(f)$ of  $\cL^2$. 
 \item[(iv)] For any open set $U\subset [0,\infty)$, the  submodules    closure of  $C_0^{\infty}(U,\Sigma\cE(\sr0))$ in $ C_0^{\infty}(U,L^2)$ are $\nt$-equivariant and  define a subsheaf $\overline{\Sigma\cE}\subset \cL^2$ of closed $\cO$-modules on $\scal2$.
 \item[(v)] For $U\subset [0,\infty)$ open, a section $\xi\in C^\infty_0(U, L^2)$ belongs to the submodule 	$C_0^{\infty}(U,\overline{\Sigma\cE})$ if and only if each Fourier component $\fourier(\xi)(n)$ as in \eqref{pointwisefourier} belongs to the closure in $C_0^{\infty}(U,\C)$ of the submodule generated by the multiples of the function $\zeta\left(\frac 12-\frac{2\pi i n}{L}\right)$  ($\zeta(z)$ is the Riemann zeta function).
 \end{enumerate}
 \end{proposition}
 
 \proof (i)~holds since the Fourier transform is  a unitary isomorphism. 

 (ii)~The  sheaf  $\cL^2$ on $[0,\infty)$  is defined by associating to an open subset $U\subset [0,\infty)$ the space  $\cF(U)=C_0^{\infty}(U,L^2)$ of smooth sections vanishing at $L=0$ of the vector bundle $L^2$.  The action of $\nt$ on $\cL^2$ is given, for $n\in \N^\times$ and for any pair of opens $U$ and $U'$ of $[0,
 \infty)$, with $n U\subset U'$, by 
 \begin{equation}\label{F(omega,n)}
\cF(U,n): C_0^{\infty}(U',L^2)\to C_0^{\infty}(U,L^2), \qquad 
\cF(U,n)(\xi)(x)= \sigma_n(\xi(nx)).
 \end{equation}
 Note that with $\mu=\exp x$ one has $\xi(nx)\in L^2(C_{\mu^n})$ and $\sigma_n(\xi(nx))\in L^2(C_\mu)$.
 By construction one has:  $\sigma_n\sigma_m=\sigma_{nm}$,  thus the above action of $\nt$ turns  $\cL^2$  into a sheaf on $\scal2=\rnt$. \newline
 (iii)~ By Lemma \ref{fouriertruncated1} (i), $\cE(f)(u)$ is pointwise well-defined, it is 
$O(u^{1/2})$ for $u\to 0$, and of rapid decay for $u\to\infty$. By (ii) of the same lemma one has 
\[
\fourier(\Sigma_\mu(\cE(f)))(n)=L^{-\frac 12}\int_{C_\mu} \Sigma_\mu(\cE(f))(u)u^{-\frac{2\pi i n}{L}}d^*u =L^{-\frac 12}\int_{\R_+^*} \cE(f)(u)u^{-\frac{2\pi i n}{L}}d^*u. 
\]
It then follows from \cite{ccspectral} (see $(6.4)$  which is valid for $z=\frac{2\pi  n}{L}\in \R$) that  
\begin{equation}\label{pointwisefourier1}
\fourier(\Sigma_\mu(\cE(f)))(n)=L^{-\frac 12}\zeta\left(\frac 12-\frac{2\pi i n}{L}\right)\int_{\R_+^*}u^{\frac 12} f(u)u^{-\frac{2\pi i n}{L}}d^*u.
\end{equation}
Since $f\in \sr0$, with $w(f)(u):=u^{1/2}f(u)$, the multiplicative Fourier transform $\fourier(w(f))=\psi$,
$\psi(z):=\int_{\R_+^*}f(u)u^{\frac 12-iz}d^*u$ is holomorphic in the complex half-plane defined by $\Im(z)>-5/2$ \cite{ccspectral}. Moreover, by construction $\sr0$ is stable under the operation $f\mapsto u\partial_u f+\frac 12 f$,  hence $w(\sr0)$ is stable under $f\mapsto u\partial_uf$. This operation  multiplies $\fourier(w(f))(z)=\psi(z)$ by $iz$. This argument shows that for any integer $m>0$, $z^m\psi(z)$ is bounded in a strip around the real axis and hence that  the  derivative $\psi^{(k)}(s)$ is $O(\vert s\vert^{-m})$ on $\R$, for any $k\geq 0$. By applying classical estimates  due to Lindelof \cite{lindelof}, (see \cite{backlund} inequality (56)), the derivatives  $\zeta^{(m)}(\frac 12+iz)$ are $O(\vert z\vert^ \alpha)$ for any $\alpha>1/4$. Thus all derivatives $\partial_L^m$ of the function \eqref{pointwisefourier1}, now re-written as $h(L,n):=L^{-\frac 12}\zeta\left(\frac 12-\frac{2\pi i n}{L}\right)\psi(\frac{2\pi n}{L})$, are sequences of rapid decay as functions of $n\in \Z$.   It follows that $\Sigma\cE(f)$ is a smooth (global) section of the vector bundle $L^2$ over $(0,\infty)$. Moreover, when $n\neq 0$ the function $h(L,n)$ tends to $0$ when $L\to 0$ and the same holds for all derivatives $\partial_L^m h(L,n)$. In fact, for any $m, k\geq 0$,  one has 
\[
\sum_{n\neq 0} \vert\partial_L^m h(L,n)\vert^2=O(L^k) \ \ \rm{when} \ \ L\to 0.
\]
This result is a consequence of  the rapid decay at $\infty$ of the derivatives of the function $\psi$, and the above estimate of $\zeta(z)$ and its derivatives. For $n=0$ one has $h(L,0)=L^{-\frac 12}\zeta(\frac 12)\psi(0)$.
\newline
 (iv)~For any open subset $U\subset [0,\infty)$ the vector space $C_0^{\infty}(U,L^2)$ admits a natural Frechet topology with generating seminorms of the form ($K\subset U$ compact subset)
 \begin{equation}\label{seminorms}
 p_K^{(n,m)}(\xi):=\max_K L^{-n}\Vert \partial_L^m \xi(L)\Vert_{L^2}.
 \end{equation}
 One obtains a space of smooth sections $C_0^{\infty}(U,\Sigma\cE)\subset C_0^{\infty}(U,L^2)$ defined as sums of products  $\sum h_j\Sigma\cE(f_j)$, with $f_j\in  \sr0$ and $h_j\in C_0^{\infty}(U,L^2)$. The map $\sigma_n:L^2(C_{\mu^n})\to L^2(C_\mu)$ in (ii) is continuous, and  from the equality
 $\sigma_n\circ \Sigma_{\mu^n}=\Sigma_\mu$ it follows (here we use the notations as in the proof of (ii)) that the sections $\xi \in C^{\infty}(U',L^2(C'))$ which belong to $C^{\infty}(U',\overline{\Sigma\cE(\sr0)})$ are mapped by $\cF(U,n)$ inside $C^{\infty}(U,\overline{\Sigma\cE(\sr0)})$. In this way one obtains a sheaf $\overline{\Sigma\cE}\subset \cL^2$ of $\cO$-modules over $\scal2$.\newline
 (v)~Let $\xi \in H^0(U,\overline{\Sigma\cE})$. By hypothesis, $\xi$ is in the closure of  $C_0^{\infty}(U,\Sigma\cE(\sr0))\subset C_0^{\infty}(U,L^2)$ for the Frechet topology. The Fourier components of $\xi$ define continuous maps in the Frechet topology, thus it follows from \eqref{pointwisefourier1} that the functions  $f_n=\fourier(\xi)(n)$  are in the closure,  for the Frechet topology on $C_0^{\infty}(U,\C)$, of $C_0^{\infty}(U,\C)g_n$, where   $g_n(L):=\zeta\left(\frac 12-\frac{ 2\pi i n }{L}\right)$ is a  multiplier of $C_0^{\infty}(U,\C)$. This conclusion holds thanks to the moderate growth of the Riemann zeta function and its derivatives on the critical line. Conversely, let $\xi\in C_0^{\infty}(U,L^2)$ be such that each of its Fourier components $\fourier(\xi)(n)$ belongs to the closure  for the Frechet topology of $C_0^{\infty}(U,\C)$, of $C_0^{\infty}(U,\C)g_n$. Let $\rho\in C_c^{\infty}([0,\infty),[0,1])$ defined to be identically equal to $1$ on $[0,1]$ and with support inside $[0,2]$. The functions $\alpha_k(x):=\rho((kx)^{-1})$ ($k>1$) fulfill the following three properties
 \begin{enumerate}
 \item 	$\alpha_k(x)=0, \quad \forall x<(2k)^{-1}$,\qquad  $\alpha_k(x)=1, \quad \forall x>k^{-1}$.
 \item $\alpha_k\in C_0^{\infty}([0,\infty))$.
 \item For all $m >0$ there exists $C_m<\infty$ such that $\vert x^{2m}\partial_x^m\alpha_k(x)\vert\leq C_m k^{-1}$ $\forall x\in [0,1], k>1$.
 \end{enumerate}
 To justify  (3), note that $x^2\partial_x f((kx)^{-1})=-k^{-1}f'((kx)^{-1})$ and that the derivatives of $\rho$ are bounded. Thus one has 
 $\vert(x^2\partial_x)^m\alpha_k(x)\vert\leq \Vert \rho^{(m)}\Vert_{\infty} k^{-m}$ $\forall x\in [0,\infty),~ k>1$ which implies (3) by induction on $m$.\newline
 Thus, when $k\to \infty$ one has $\alpha_k\xi \to \xi$ in the Frechet topology of $C_0^{\infty}(U,L^2)$. This is clear if $0\notin U$ since then, on any  compact subset $K\subset U$, all $\alpha_k$ are identically equal to $1$ for $k>(\min K)^{-1}$. Assume now that $0\in U$ and let $K=[0,\epsilon]\subset U$. With the notation of \eqref{seminorms} let us show that $p_K^{(n,m)}((\alpha_k-1)\xi)\to 0$ when $k\to \infty$. Since $\alpha_k(x)=1, \, \forall x>k^{-1}$ one has, using the finiteness of $p_K^{(n+1,m)}(\xi)$
 \[
  \max_K L^{-n} \Vert ((\alpha_k-1) \partial_L^m \xi)(L)\Vert_{L^2}\leq 
  \max_{(0,1/k]} L^{-n} \Vert (\partial_L^m \xi)(L)\Vert_{L^2}\stackrel{k\to \infty}{\to} 0.
 \]
 Then one obtains
\[
L^{-n} \partial_L^m((\alpha_k-1) \xi)(L)=L^{-n}((\alpha_k-1) \partial_L^m \xi)(L)+\sum_1^m L^{-n} {m\choose j}(\partial_L^j \alpha_k)(L) (\partial_L^{m-j} \xi)(L).
\]
 Thus using (3) above and the finiteness of the norms $p_K^{(n+2j,m-j)}(\xi)$ one derives: $p_K^{(n,m)}((\alpha_k-1)\xi)\to 0$ when $k\to \infty$. It remains to show that $\alpha_k \xi$  belongs to the submodule 	$C_0^{\infty}(U,\overline{\Sigma\cE})$. It is enough to show that for $K\subset (0,\infty)$  a compact subset with $\min K>0$, one can approximate $\xi$  by elements of $C_0^{\infty}(U,\Sigma\cE)$ for the norm $p_K^{(0,m)}$. Let $P_N$ be the orthogonal projection in $L^2(C_\mu)$ on the finite-dimensional subspace determined by the vanishing of all Fourier components $\fourier(\xi)(\ell)$ for any $\ell, \vert \ell\vert>N$. Given $L\in K$ and $\epsilon >0$ there exists $N(L,\epsilon)<\infty$ such that 
 \begin{equation}\label{approxpn}
 \Vert (1-P_N)\partial_L^j\xi(L)\Vert <\epsilon \quad \forall j\leq m, \ N\geq N(L,\epsilon).
 \end{equation}
 The smoothness of $\xi$ implies that there exists an open neighborhood $V(L,\epsilon)$ of $L$ such that \eqref{approxpn} holds in $V(L,\epsilon)$. The compactness of $K$ then shows that there exists a finite $N_K$ such that 
 \begin{equation*}\label{approxpn1}
 \Vert (1-P_N)\partial_L^j\xi(L)\Vert <\epsilon \quad \forall j\leq m, \ L\in K, \ N\geq N_K.
 \end{equation*} 
 It now suffices to show that one can approximate $P_N\xi$, for the norm $p_K^{(0,m)}$, by elements of $C_0^{\infty}(U,\Sigma\cE)$. To achieve this result, we let $L_0\in K$ and $\delta_j\in C_c^{\infty}(\R_+^*)$, $\vert j \vert \leq N$  be such that 
 \[
 \fourier(\delta_j)\left(\frac{ 2 \pi j'}{L_0}\right)=\begin{cases}1 & j=j'\\
 0 & j\neq j'\end{cases}\quad \forall j', \vert j' \vert \leq N, \ \ \int_{\R_+^*}u^{1/2}\delta_j(u)d^*u=0\quad \forall j, \vert j \vert \leq N.
 \]
 One construct $\delta_j$  starting with a function $h\in C_c^{\infty}(\R_+^*)$ such that
 $\fourier(h)\left(\frac{ 2 \pi j}{L_0}\right)\neq 0$  and acting on $h$ by a differential polynomial whose effect is to multiply $\fourier(h)$ by a polynomial vanishing on all $\frac{ 2 \pi j'}{L_0}$, $j'\neq j$ and at $i/2$. By hypothesis each Fourier component $\fourier(\xi)(n)$ belongs to the closure in $C_0^{\infty}(U,\C)$ of the multiples of the function $\zeta\left(\frac 12-\frac{2\pi i n}{L}\right)$. Thus, given $\epsilon>0$ one has functions $f_n\in C_0^{\infty}(U,\C)$, $\vert n\vert \leq N$ such that \[
 \max_K \vert\partial_L^j\left(\fourier(\xi)(n)-\zeta\left(\frac 12-\frac{2\pi i n}{L}\right)f_n(L)\right)\vert \leq \epsilon\quad \forall j\leq m, \quad \vert n\vert \leq N.
 \]
 We now can find a small open neighborhood $V$ of $L_0$ and functions $\phi_j\in C^{\infty}(V)$, $\vert j\vert \leq N$ such that 
 \begin{equation}\label{approxpn2}
 \sum \phi_j(L)\fourier(\delta_j)\left(\frac{ 2 \pi n}{L}\right)=L^{\frac 12}f_n(L)\quad \forall L\in V.
 \end{equation}
 This is possible because the determinant of the matrix $M_{n,j}(L)=\fourier(\delta_j)\left(\frac{ 2 \pi n}{L}\right)$ is non-zero in a neighborhood of $L_0$ where $M_{n,j}(L_0)$ is the identity matrix. The even functions $d_j(u)$ on $\R$, which agree with $u^{-1/2}\delta_j(u)$ for $u>0$, are all in $\sr0$ since $\int_\R d_j(x)dx=2\int_{\R_+^*}u^{1/2}\delta_j(u)d^*u=0$. One then has \[
\fourier(\Sigma_\mu(\cE(d_j)))(n)=L^{-\frac 12}\zeta\left(\frac 12-\frac{2\pi i n}{L}\right)\fourier_\mu(\delta_j)\left(\frac{ 2 \pi n}{L}\right) 
\]
 by \eqref{pointwisefourier1}, and by \eqref{approxpn2} one  gets
 \[
 \sum \phi_j(L)\fourier(\Sigma_\mu(\cE(d_j)))(n)=\zeta\left(\frac 12-\frac{2\pi i n}{L}\right)f_n(L)\quad \forall L\in V.
\]
 One finally covers $K$ by finitely many such open sets $V$ and use a partition of unity subordinated to this covering to obtain smooth functions $\varphi_\ell\in C_c^{\infty}(0,\infty)$, $g_\ell \in \sr0$ such that the Fourier component of index $n$, $\vert n\vert \leq N$, of $\sum  \varphi_\ell \Sigma\cE(g_\ell)$ is  equal  to  $\zeta\left(\frac 12-\frac{2\pi i n}{L}\right)f_n(L)$ on $K$. This shows that $\xi$ belongs to the closure of  $C_0^{\infty}(U,\Sigma\cE(\sr0))\subset C_0^{\infty}(U,L^2)$. 
 \endproof

 We recall that the space of  global sections $H^0(\cT,\cF)$ of a sheaf of sets $\cF$ in a Grothendieck topos $\cT$  is defined to be the set $\Hom_{\cT}(1,\cF)$, where $1$ denotes the terminal object of $\cT$.  For $\cT=\scal 2$ and $\cF$  a   sheaf of sets on $[0,\infty)$,  $1$ assigns to an open set $U\subset [0,\infty)$ the single element $*$, on which $\nt$ acts as the identity. Thus, we understand an element of $\Hom_{\scal2}(1,\cF)$ as  a global section $\xi$ of $\cF$, where $\cF$ is viewed as a sheaf on $[0,\infty)$  invariant under the action of $\nt$.\vspace{.05in} 
 
With the notations of Proposition~\ref{sheaves} and  for  $\xi\in\Hom_{\scal2}(1,\cL^2)$, we  write $\widehat \xi(L,n):=\fourier(\xi)(n)$ for the (multiplicative) Fourier components of $\xi$. Then we have

\begin{lemma}\label{sheaves1} The following facts hold
 \begin{enumerate}
     \item[(i)] The map 
  \begin{equation}\label{mapgamma}
 \gamma:H^0(\scal2,\cL^2)\to C_0^{\infty}([0,\infty),\C)\times C_0^{\infty}([0,\infty),\C)\qquad \gamma(\xi)=\widehat \xi(\pm 1)
\end{equation}
is an isomorphism of $\C$-vector spaces.
 \item[(ii)] The subspace $\gamma(H^0(\scal2,\overline{\Sigma\cE}))$ is the closed ideal (for the Frechet topology on $C_0^{\infty}([0,\infty),\C)$), generated by the multiplication with the functions $\zeta\left(\frac 12 \mp\frac{ 2\pi i }{L}\right)$.
 \end{enumerate}
 \end{lemma}
\proof (i)~Let $\xi\in\Hom_{\scal2}(1,\cL^2)$: this is a global section $\xi \in C_0^{\infty}([0,\infty),L^2)$ invariant under the action of $\nt$, \ie such that  $\sigma_n(\xi(nL))=\xi(L)$ for all pairs $(L,n)$. The Fourier components $\widehat \xi(L,n)$ of any such section are smooth functions of $L\in [0,\infty)$ vanishing  at $L=0$, for $n\neq 0$, as well as all their derivatives. The equality $\sigma_n(\xi(L))=\xi(L/n)$ entails, for $n>0$, 
\begin{align*}
\widehat \xi(L,n)&=L^{-\frac 12}\int_{C_\mu} \xi(L)(u)u^{-\frac{2\pi i n}{L}}d^*u=L^{-\frac 12}\int_{C_{\mu}^{1/n}} \xi(L/n)(u)u^{-\frac{2\pi i }{L}}d^*u=\\&=n^{-\frac 12}\widehat \xi(L/n,1).
\end{align*}
 This shows that the $\widehat \xi(L,n)$ are uniquely determined, for $n>0$ by the function $\widehat \xi(L,1)$ and,  for $n<0$, by the function $\widehat \xi(L,-1)$. With $g(L)=\widehat \xi(L,0)$ one has: $g(L)=n^{-\frac 12}g(L/n)$ for all $n>0$. This implies, since $\Q_+^*$ is dense in $\R_+^*$ and $g$ is assumed to be smooth, that $g$ is proportional to $L^{-\frac 12}$ and hence  identically $0$, since it corresponds to a global section smooth at $0\in [0,\infty)$. This argument proves that  $ \gamma$ is injective. Let us show that $\gamma$ is also surjective. Given a pair of functions $f_\pm\in C_0^{\infty}([0,\infty),\C)$ we  construct a global section $\xi\in H^0(\scal2,\cL^2)$ such that $\gamma(\xi)=(f_+,f_-)$. One defines $\xi(L)\in L^2(C_\mu)$ by by means of its Fourier components  set to be $ \widehat \xi(L,0):=0$, and  for $n\neq 0$ by
 \[
 \widehat \xi(L,n):=\vert n\vert^{-\frac 12}f_{\rm sign(n)}(L/n).
\]
 Since $f_\pm(x)$ are of rapid decay for $x\to 0$,  $\sum \vert \widehat \xi(L,n)\vert^2<\infty$, thus $\xi(L)\in L^2(C_\mu)$. All derivatives of $f_\pm(x)$ are also of rapid decay for $x\to 0$, thus all derivatives $\partial_L^k(\xi(L))$ belong to $L^2(C_\mu)$ and that the $L^2$-norms $\Vert\partial_L^k(\xi(L))\Vert$ are of rapid decay for $L\to 0$.  By construction  $\sigma_n(\xi(L))=\xi(L/n)$, which entails  $\xi\in H^0(\scal2,\cL^2)$ with $\gamma(\xi)=(f_+,f_-)$.\newline
 (ii)~Let $\xi \in H^0(\scal2,\overline{\Sigma\cE})$. By Proposition \ref{sheaves} (v),   the functions  $f_\pm=\widehat \xi(L,\pm 1)$ are in the closure,  for the Frechet topology on $C_0^{\infty}([0,\infty),\C)$, of the ideal generated by the functions $\zeta\left(\frac 12 \mp\frac{ 2\pi i }{L}\right)$. Conversely, let $\xi \in H^0(\scal2,\cL^2)$ and assume that $\gamma(\xi)$ is in the closed submodule generated by multiplication with $\zeta\left(\frac 12 \mp \frac{ 2\pi i}{L}\right)$. The $\nt$-invariance of $\xi$ implies $\widehat \xi(L,n)=\vert n\vert^{-\frac 12}\widehat \xi(L/\vert n\vert,{\rm sign}(n))$ for $n\neq 0$. Thus the Fourier components $\widehat \xi(L,n)$ belong to the closure in $C_0^{\infty}(U,\C)$ of the multiples of the function $\zeta\left(\frac 12-\frac{2\pi i n}{L}\right)$, then Proposition \ref{sheaves} (v) again implies $\xi \in H^0(\scal2,\overline{\Sigma\cE})$.\endproof 
 
 The action of  $\R_+^*$  on the  sheaf $\cL^2$  is given by the action  $\rep$  on the Fourier components of its  sections $\xi$. With $\mu=\exp L$, $L\in (0,\infty)$, $n\in\N^*$ and $\lambda\in\R^*_+$, this is 
 \begin{align}\label{action}
 \widehat {\rep(\lambda)\xi}(L,n)&=L^{-\frac 12}\int_{C_\mu} \xi(\lambda^{-1}u)u^{-\frac{2\pi i n}{L}}d^*u =\lambda^{-\frac{2\pi i n}{L}} L^{-\frac 12}\int_{C_\mu} \xi(v)v^{-\frac{2\pi i n}{L}}d^*v=\\&=\lambda^{-\frac{2\pi i n}{L}}\widehat \xi(L,n)\notag
 \end{align}
 
The following result  explains in particular how the quotient sheaf $\cL^2/\overline{\Sigma\cE}$ on $\scal2$  handles eventual multiplicities of critical zeros of the zeta function.
   
 \begin{theorem}
 The  induced action of $\R_+^*$ on the global sections $H^0(\scal2,\cL^2/\overline{\Sigma\cE})$   is canonically isomorphic to the  action of $\R_+^*$,  via multiplication with $\lambda^{is}$, on the quotient of the Schwartz space $\cS(\R)$ by the closure of the ideal generated by multiples of $\zeta\left(\frac 12 +is\right)$. 
 \end{theorem}
 \proof 
 We first show that the canonical map $q:H^0(\scal2,\cL^2)\to H^0(\scal2,\cL^2/\overline{\Sigma\cE}))$ is surjective. Let $\xi \in H^0(\scal2,\cL^2/\overline{\Sigma\cE}))$: as  a section of $\cL^2/\overline{\Sigma\cE}$ on $[0,\infty)$,  there exists an open neighborhood $V=[0,\epsilon)$ of $0\in [0,\infty)$ and a section $\eta\in C_0^{\infty}(V,L^2)$ such that the class of $\eta$ in $C_0^{\infty}(V,L^2/\overline{\Sigma\cE})$  is the restriction of $\xi$ to $V$. The Fourier components $\widehat \eta(L,n)$ are meaningful for $L\in V$. Since $\xi$ is $\nt$-invariant, for any $n\in \nt$ the class of $\cF(V/n,n)(\eta)$, with $\cF(V/n,n)(\eta)(L):= \sigma_n(\eta(nL))$ (see \eqref{sigman}) is equal to the class of the restriction of $\eta$ in $C_0^{\infty}(V/n,L^2/\overline{\Sigma\cE})$. We thus obtain 
\[
 \eta(L) -\cF(V/n,n)(\eta)\in C_0^{\infty}(V/n,\overline{\Sigma\cE})
 \]
 Furthermore, the Fourier components of $\alpha=\cF(V/n,n)(\eta)$ are given by 
 \[
 \widehat \alpha(L,k)=n^{\frac 12}\widehat \eta(nL,nk).
 \]
 Thus, the  functions 
 \[
 \widehat \eta(L,k)-\widehat \alpha(L,k)=  \widehat \eta(L,k)-n^{\frac 12}\widehat \eta(nL,nk)
 \]
 are the Fourier components of an element in  $C_0^{\infty}(V/n,\overline{\Sigma\cE})$. 
  By Proposition \ref{sheaves} (v),  these components are in the closure  in $C_0^{\infty}(V/n,\C)$ of the multiples of the function $\zeta\left(\frac 12-\frac{2\pi i k}{L}\right)$. For $k=1$  the function $\widehat \eta(L,1)-n^{\frac 12}\widehat \eta(nL,n)$ is in the closure  in $C_0^{\infty}(V/n,\C)$ of the multiples of the function $\zeta\left(\frac 12-\frac{2\pi i }{L}\right)$. This implies that  $\widehat \eta(L,n)-n^{-\frac 12}\widehat \eta(L/n,1)$ is in the closure  in $C_0^{\infty}(V,\C)$ of the multiples of the function $\zeta\left(\frac 12-\frac{2\pi in }{L}\right)$. For $k=-1$ one obtains  a similar result for the Fourier components $\widehat \eta(L,n)$, $n<0$. Thus, again by Proposition \ref{sheaves} (v), one knows that  the class of $\eta$  in $C_0^{\infty}(V,L^2/\overline{\Sigma\cE}))$ s not altered if one replaces $\eta$ with $\eta_1\in C_0^{\infty}(V,L^2)$  
  \[
  \widehat \eta_1(L,n):=\vert n\vert^{-\frac 12}\eta(L/\vert n\vert,{\rm sign}(n)).
  \]
 Next step is to extend the functions $\eta(L,\pm 1)\in C_0^{\infty}(V,\C)$ to  $f_\pm\in  C_0^{\infty}([0,\infty),\C)$  fulfilling the following property. For any open set $U\subset [0,\infty)$ and a section $\beta\in C_0^{\infty}(U,L^2)$, with the class of $\beta$ in $C_0^{\infty}(U,L^2/\overline{\Sigma\cE}))$ being the restriction of $\xi$ to $U$, the functions $
 \widehat \beta(L,\pm 1)-f_\pm(L)
 $ 
 belong to the closure  in $C_0^{\infty}(U,\C)$ of the multiples of the function $\zeta\left(\frac 12\mp\frac{2\pi i }{L}\right)$. To construct $f_\pm$ one considers the sheaf $\cG_\pm$ which is the quotient of the sheaf of $C_0^{\infty}([0,\infty),\C)$ functions by the closure of the ideal subsheaf generated by the multiples of the function $\zeta\left(\frac 12\mp\frac{2\pi i }{L}\right)$. Since the latter is a module over the sheaf of $C^{\infty}$ functions, it is a fine sheaf, thus a global section  of $\cG_\pm$ can be lifted to a function. By Proposition \ref{sheaves} (v), the Fourier components $\widehat \xi_j(L,\pm 1)$ of local sections $\xi_j$ of $L^2$ representing $\xi$ define a global section of $\cG_\pm$. The functions $f_\pm$ are obtained by lifting these sections. By appealing to Lemma \ref{sheaves1}, we let $\phi \in H^0(\scal2,L^2)$ to be the unique global section such that $\gamma(\phi)=(f_+,f_-)$. Then we show that $q(\phi)=\xi$.  We have already proven that the restrictions to $V=[0,\epsilon)$ are the same. Thus it is enough to show that given $L_0>0$ and a lift $\xi_0\in C_0^{\infty}(U,L^2)$ of $\xi$ in a small open interval $U$ containing $L_0$, the difference $\delta=\phi-\xi_0$ is a section of $\overline{\Sigma\cE}$. Again by Proposition \ref{sheaves} (v), it suffices to show that the Fourier components  $\widehat \delta(L,n)$ are in the closure of the ideal generated by multiples of $\zeta\left(\frac 12-\frac{2\pi in }{L}\right)$. The $\nt$-invariance of $\xi$ shows that $\cF(U/n,n)(\xi_0)$ (see \eqref{F(omega,n)}) is a lift of $\xi$ in $U/n$. Thus by the defining properties of the functions $f_\pm$ one has 
 \[
 \widehat{\cF(U/n,n)(\xi_0)}(\pm 1)-f_\pm \in \overline{C^{\infty}(U,\C)\zeta_\pm}, \qquad\text{for}\quad   \zeta_\pm(L)=\zeta\left(\frac 12\mp\frac{2\pi i }{L}\right).
 \]
  With a similar argument and using the invariance of $\phi$ under the action of $\cF(U/n,n)$,    one obtains that $\widehat \delta(n)$ is in the closure of the ideal generated by the   multiples of $\zeta\left(\frac 12-\frac{2\pi in }{L}\right)$.\newline
  We have thus proved that $q:H^0(\scal2,\cL^2)\to H^0(\scal2,\cL^2/\overline{\Sigma\cE}))$ is surjective. One then obtains the exact sequence 
\begin{equation}\label{exactseq}
  0\to H^0(\scal2,\overline{\Sigma\cE}))\to H^0(\scal2,\cL^2)\to H^0(\scal2,\cL^2/\overline{\Sigma\cE}))\to 0.
  \end{equation}
 This sequence is equivariant for the action \eqref{action} of $\rep$ of $\R_+^*$ on the bundle $L^2$. 
 For $h\in L^1(\R_+^*,d^*u)$ one has 
\begin{equation}\label{equivariant}
  \widehat {(\rep(h)\xi)}(L,n)=\fourier(h)\left(\frac{2\pi  n}{L}\right) \widehat \xi(L,n).
 \end{equation}
To obtain the required isomorphism between the two spectral realizations, one uses  the isomorphism \eqref{mapgamma} of Lemma \ref{sheaves1}. Denote for short $(C_0^{\infty})^2=C_0^{\infty}([0,\infty))\times C_0^{\infty}([0,\infty))$.
 One maps the Schwartz space $\cS(\R)$ to $(C_0^{\infty})^2$ by 
 \begin{equation*}\label{mapschwartz}
  \Phi(f)=(\Phi_+(f),\Phi_-(f)), \qquad \Phi_\pm(f)(L):=f\left(\pm \frac{2\pi}{L}\right).
    \end{equation*}
   $ \Phi$ is well defined since  all derivatives of $\Phi_\pm(f)(L)$ tend to $0$ when $L\to 0$ (any function $f\in\cS(\R)$ is of rapid decay as well as all its derivatives). 
 The exact sequence \eqref{exactseq}, together with Lemma \ref{sheaves1}, then gives an induced isomorphism  
 \begin{equation*}\label{mapgamma1}
 \gamma: H^0(\scal2,\cL^2/\overline{\Sigma\cE}))\simeq (C_0^{\infty})^2/ \left(\overline{C_0^{\infty}\zeta_+}\times \overline{C_0^{\infty}\zeta_-}\right).
    \end{equation*}
 In turn, the map $ \Phi$ induces a morphism 
  \begin{equation*}\label{mapPhi}
 \tilde \Phi:\cS(\R)/(\cS(\R)\zeta) \to (C_0^{\infty})^2/ \left(\overline{C_0^{\infty}\zeta_+}\times \overline{C_0^{\infty}\zeta_-}\right).
  \end{equation*}
  By \eqref{equivariant}  this morphism is equivariant for the action of $\R_+^*$. The map $ \Phi$ is not an isomorphism since elements of its range have finite limits at $\infty$. However it is injective and its range contains all elements of  $(C_0^{\infty})^2$ which have compact support.  Since $\zeta_\pm(L)=\zeta\left(\frac 12\mp\frac{2\pi i }{L}\right)$ tends to a finite non-zero limit when $L\to 0$,   $\tilde \Phi$ is an isomorphism.\endproof 
 
 \begin{remark}\label{remtauto}  By a Theorem of Whitney (see \cite{Malgrange}, Corollary 1.7), the closure  of the ideal of multiples of $\zeta\left(\frac 12 +is\right)$ in $\cS(\R)$ is the subspace of those  $f\in\cS(\R)$ which vanish of the same order as $\zeta$ at every (critical) zero $s\in Z$. Thus if any such zero is a multiple zero of order $m>1$, one finds that the action of $\R_+^*$ on the global sections of the quotient sheaf $\cL^2/\overline{\Sigma\cE}$ admits a non-trivial Jordan decomposition  of the form
 $$
 \rep(\lambda)\xi=\lambda^{is}(\xi+N(\lambda)\xi), \ \ 
 $$
 with $N(\lambda)^m=0$ and  $(1+N(u))(1+N(v))=1+N(uv)$ for all $u,v\in \R_+^*$.
 \end{remark}

\end{document}